\documentstyle[11pt]{article}

\newcommand{\bbq}{{\mathbf Q}}

\newcommand{\qed}{\hfill$\Box$}
\newcommand{\inv}{{\mathrm inv}}
\newcommand{\Des}{{\mathrm Des}}
 
\newcommand{\maj}{{\mathrm maj}}
\newcommand{\ldes}{{\mathrm ldes}}
\newcommand{\lind}{{\mathrm lind}}
\newcommand{\Peak}{{\mathrm Peak}}
\newcommand{\tail}{{\mathrm tail}}
\newcommand{\bpi}{{\bar{\pi}}}
\newcommand{\hpi}{{\hat{\pi}}}
\newcommand{\bt}{{\bar{t}}}
\newcommand{\hf}{{\hat{f}}}
\newcommand{\sign}{{\mathrm sign}}
\newcommand{\SYT}{{\cal SYT}}
\newcommand{\DP}{{\cal P}}

\newtheorem{thm}{Theorem}[section]

\newtheorem{lem}[thm]{Lemma}

\newtheorem{cor}[thm]{Corollary}

\newtheorem{prb}[thm]{Problem}

\begin{document}

\title{Equidistribution and Sign-Balance on $321$-Avoiding Permutations}
\bibliographystyle{acm}
\author{Ron M.\ Adin%
\thanks{Department of Mathematics and Statistics, Bar-Ilan University,
Ramat-Gan 52900, Israel. 
Email: {\tt radin@math.biu.ac.il} } 
\and Yuval Roichman%
\thanks{Department of Mathematics and Statistics, Bar-Ilan University,
Ramat-Gan 52900, Israel. 
Email: {\tt yuvalr@math.biu.ac.il} } 
\thanks{Both authors partially supported by the EC's IHRP Programme, 
within the Research Training Network ``Algebraic Combinatorics in Europe'', 
grant HPRN-CT-2001-00272.}}
\date{Submitted: May 5, 2003; Revised: January 12, 2004}

\maketitle

\begin{abstract}
Let $T_n$ be the set of $321$-avoiding permutations of order $n$.
Two properties of $T_n$ are proved:
(1) The {\em last descent} and {\em last index minus one} statistics are 
equidistributed over $T_n$, and also over subsets of permutations 
whose inverse has an (almost) prescribed descent set.
An analogous result holds for Dyck paths.
(2) The sign-and-last-descent enumerators for $T_{2n}$ and $T_{2n+1}$ are 
essentially equal to the last-descent enumerator for $T_n$.
The proofs use a recursion formula for an appropriate multivariate generating 
function.
\end{abstract}

\section{Introduction}

\subsection{Equidistribution}

One of the frequent themes in combinatorics is identifying two 
distinct parameters on the same set which are {\it equidistributed},
i.e., share the same generating function.
The first substantial result of this kind on permutations,
by MacMahon~\cite{MM}, received a remarkable refinement by Foata
and Sch\"utzenberger~\cite{FS} (see also~\cite{G, GG}).
They proved the equidistribution of inversion number and major index,
not only on the whole symmetric group, but also on distinguished subsets
of permutations (those whose inverses have a prescribed descent set). 

Equidistribution theorems on descent classes were shown to be
closely related to the study of polynomial rings; 
see, e.g., \cite{R, ABR}.
Motivated by the properties of certain quotient rings,
studied by Aval and Bergeron~\cite{AB, ABB} (see also Section~\ref{remarks}
below), we looked for an analogue for the set of $321$-avoiding permutations
of the above mentioned theorem of Foata and Sch\"utzenberger.

\medskip

Let $S_n$ be the symmetric group on $n$ letters, and let
$$
T_n := \{\pi\in S_n \,|\, \not\exists\, i<j<k \mbox{\rm\ such that\ } \pi(i) > \pi(j) > \pi(k)\}
$$
be the set of $321$-avoiding (or two-row shaped) permutations in $S_n$.
For $\pi \in T_n$ define the following statistics:
\begin{eqnarray*}
\inv(\pi) 
&:=& \mbox{\rm inversion number of $\pi$}\\ 
&(=& \mbox{\rm length of $\pi$ w.r.t. the usual generators of $S_n$})\\
\ldes(\pi) 
&:=& \mbox{\rm last descent of $\pi$} 
     = \max\{1\le i\le n-1 \,|\, \pi(i) > \pi(i+1)\}\\
& & (\mbox{\rm where $\ldes(\pi) := 0$ for $\pi=id$})\\
\lind(\pi) 
&:=& \pi^{-1}(n), 
     \mbox{\rm\ the index of the digit ``$n$'' in $\pi$}\\
\Des(\pi^{-1}) 
&:=& \mbox{\rm descent set of\ } \pi^{-1} 
     = \{1\le i\le n-1 \,|\, \pi^{-1}(i) > \pi^{-1}(i+1)\}
\end{eqnarray*}

The following theorem is a $T_n$-analogue of~\cite[Theorem 1]{FS}.

\begin{thm}\label{t.equi}
The statistics `` $\ldes$'' and `` $\lind-1$'' are equidistributed over $T_n$.
Moreover, for any $B\subseteq[n-2]$ they are equidistributed over 
the set 
$$
T_n(B):=\{\pi\in T_n \,|\, \Des(\pi^{-1})\cap [n-2] = B\},
$$
namely:
$$
\sum_{\pi\in T_n(B)} q^{\ldes(\pi)} = 
\sum_{\pi\in T_n(B)} q^{\lind(\pi)-1}.
$$
\end{thm}

%

\bigskip

Let $\DP_n$ be the set of Dyck paths of length $2n$.
For $p = (p_1,\ldots,p_{2n})\in \DP_n$, let 
\begin{eqnarray*}
\ldes(p) 
&:=& \max\{1\le i\le n-1 \,|\, p_i=+1, p_{i+1}=-1\}\\
&  &(\mbox{\rm where $\ldes(p) := 0$ for $p=(+\ldots+-\ldots-)$})\\
\lind(p) 
&:=& \max\{1\le i\le n \,|\, p_i=+1, p_{i+1}=-1\}\\
\Des(p^{-1})
&:=& \{1\le i\le n-1 \,|\, p_{2n-i}=+1, p_{2n-i+1}=-1\} 
\end{eqnarray*}

\begin{thm}\label{t.Dyck2}
$$
\sum\limits_{\pi\in T_n} q^{\ldes(\pi)}=
\sum\limits_{p\in \DP_n} q^{\ldes(p)}= 
\sum\limits_{p\in \DP_n} q^{\lind(p)-1}. 
$$
Moreover, for any $B \subseteq [n-2]$,
$$
\sum_{\pi\in T_n(B)} q^{\ldes(\pi)} =
\sum_{p\in \DP_n(B)} q^{\ldes(p)} = 
\sum_{p\in \DP_n(B)} q^{\lind(p)-1},
$$
where $T_n(B)$ is as in Theorem~\ref{t.equi} and
$$
\DP_n(B) := \{p\in \DP_n \,|\, \Des(p^{-1}) \cap [n-2] = B\}.
$$
\end{thm}

%

\subsection{Signed Enumeration}
%

Sign balance for linear extensions of posets was studied 
in~\cite{Ruskey, White, St.balance}.
The study of sign balance for pattern-avoiding permutations started with 
Simion and Schmidt~\cite{SiS}, who proved that the numbers of
even and odd $123$-avoiding permutations in $S_n$ are equal if $n$ is even,
and differ (up to a sign) by a Catalan number if $n$ is odd.
Refined sign balance (i.e., a generating function taking into account 
the signs of permutations) was investigated by 
D\'esarmenien and Foata~\cite{DF} and by Wachs~\cite{Wachs}. 
A beautiful formula for the signed Mahonian (i.e., the sign-and-major-index 
enumerator on $S_n$) was found by Simion and Gessel~\cite[Cor.~2]{Wachs}; 
see also~\cite{AGR}.

The last descent statistic is a $T_n$-analogue of the major index on 
$S_n$, as demonstrated by Theorem~\ref{t.equi} above; see also the discussion 
in Section~\ref{remarks}. Thus, a $T_n$-analogue of the signed Mahonian
on $S_n$ is the sum
$$
\sum\limits_{\pi\in T_n} \sign(\pi) q^{\ldes(\pi)}.
$$
For this sum we prove

\begin{thm}\label{t.ldes.balance}
\begin{eqnarray*}
\sum_{\pi\in T_{2n+1}}\sign(\pi) \cdot q^{\ldes(\pi)}
&=& \sum_{\pi\in T_n} q^{2\cdot \ldes(\pi)}\qquad(n\ge 0),\\
\sum_{\pi\in T_{2n}}\sign(\pi)\cdot q^{\ldes(\pi)} 
&=& (1-q)\sum_{\pi\in T_n} q^{2\cdot \ldes(\pi)}\qquad(n\ge 1).
\end{eqnarray*}
\end{thm}

This result refines (and is actually implicit in) some results in~\cite{SiS}.
Our approach is different from that of~\cite{SiS}; see Subsection~\ref{tools} 
below. After the archive posting of the current paper, 
Theorem~\ref{t.ldes.balance} has been followed by various extensions, 
analogues and refinements:
\begin{itemize}
\item[(1)] A. Reifegerste gave a bivariate refinement~\cite[Cor.~4.4]{Reife}.
\item[(2)] T. Mansour proved an analogue for $132$-avoiding permutations~\cite{Mansour}.
\item[(3)] M. Shynar found analogues for the signed enumeration of tableaux 
with respect to certain natural statistics~\cite{Shynar}.
\end{itemize}

The following phenomenon appears, with small variations, in all these results:
the standard enumerator of objects of size $n$ is essentially equal to 
the signed enumerator of objects of size $2n$. This phenomenon deserves
further study; it seems to be related to the cyclic sieving phenomenon
of~\cite{RSW}.

%
%

\subsection{Recursion}\label{tools}

Our main tool in proving the above results is a new recursive formula
for the multivariate generating function, on $321$-avoiding permutations,
involving the inversion number, the last descent, the last index, and 
the descent set of the inverse (Theorem~\ref{t.rec} below).

It should be noted that, shortly after the archive posting of this paper,
bijective proofs of Theorems~\ref{t.equi} and~\ref{t.ldes.balance} 
were presented by A.\ Reifegerste~\cite{Reife}.

\section{A Recursion Formula}\label{s.gf}

Define the multivariate generating function
\begin{equation}\label{e.gf}
f_n(\bt,x,y,z) := \sum_{\pi\in T_n} t_{\Des(\pi^{-1})}
                  x^{\inv(\pi)} y^{\ldes(\pi)} z^{\lind(\pi)},
\end{equation}
where $\bt = (t_1,t_2,\ldots)$ and $t_D := \prod_{i\in D} t_i$ for 
$D\subseteq [n-1]$.
\begin{thm}\label{t.rec}{\bf (Recursion Formula)}
$$
f_1(\bt,x,y,z) = z
$$
and, for $n\ge 2$,
\begin{eqnarray*}
(x-yz)f_n(\bt,x,y,z) 
&=& t_{n-1}x^{n}yz \cdot f_{n-1}(\bt,x,yz/x,1)\\
& & + \,(1-t_{n-1})x^{n}yz \cdot f_{n-1}(\bt,x,1,yz/x)\\
& & + \,(x-yz)z^{n} \cdot f_{n-1}(\bt,x,y,1)\\ 
& & - \,xy^{n}z^{n} \cdot f_{n-1}(\bt,x,1,1)
\end{eqnarray*}
\end{thm}

\noindent
{\bf Proof.}
The case $n=1$ is clear. Assume $n\ge 2$.
\par\noindent 
Given a permutation
$$
\pi = (\pi(1),\ldots,\pi(n-1))\in T_{n-1},
$$
insert the digit $n$ between the $k$th and $(k+1)$st places in $\pi$ 
($0\le k\le n-1$) to get a permutation
$$
\bpi = (\pi(1),\ldots,\pi(k),n,\pi(k+1),\ldots,\pi(n-1))\in S_n.
$$
Clearly, a necessary and sufficient condition for $\bpi\in T_n$ is
that $\pi$ has no descent after the $k$th place, i.e.,
\begin{equation}
\ldes(\pi)\le k\le n-1.
\end{equation}
If $\pi$ has no descents, i.e. $\pi=id$, this still holds with $\ldes(\pi):=0$.
The new statistics for $\bpi$ are:
\begin{eqnarray*}
\inv(\bpi) &=& \inv(\pi) + (n-1-k)\\
\ldes(\bpi) &=& \cases{k+1,& if $k<n-1$;\cr 
                       \ldes(\pi),& if $k=n-1$.}\\
\lind(\bpi) &=& k+1\\
\Des(\bpi^{-1}) 
&=& \cases{\Des(\pi^{-1}),& if $\pi^{-1}(n-1)\le k$;\cr 
           \Des(\pi^{-1})\cup\{n-1\},& if $k<\pi^{-1}(n-1)$.}
\end{eqnarray*}

Note that, for any $\pi\in T_{n-1}$,
$$
\ldes(\pi) \le \pi^{-1}(n-1) \le n-1;
$$
and, actually, exactly one of the following two cases holds: either
$$
\ldes(\pi) = \pi^{-1}(n-1) < n-1
\qquad(\mbox{\rm if $n-1$ is not the last digit in $\pi$})
$$
or
$$
\ldes(\pi) < \pi^{-1}(n-1) = n-1
\qquad(\mbox{\rm if $n-1$ is the last digit in $\pi$}).
$$

We can therefore compute
\begin{eqnarray*}
f_n &=& f_n(\bt,x,y,z)\\ 
&=& \sum_{\pi\in T_{n-1}} \sum_{k=\ldes(\pi)}^{n-1} 
    t_{\Des(\bpi^{-1})} x^{\inv(\bpi)} y^{\ldes(\bpi)} z^{\lind(\bpi)}\\
&=& \sum_{\pi\in T_{n-1}} t_{\Des(\pi^{-1})} x^{\inv(\pi)} 
\cdot \left[\sum_{k=\ldes(\pi)}^{\pi^{-1}(n-1)-1} 
      t_{n-1} x^{n-1-k} y^{k+1} z^{k+1}\right.\\
& & + \left.\sum_{k=\pi^{-1}(n-1)}^{n-2} 
      x^{n-1-k} y^{k+1} z^{k+1} + y^{\ldes(\pi)} z^{n}\right]\\
&=& \sum_{\pi\in T_{n-1}} t_{\Des(\pi^{-1})} x^{\inv(\pi)} 
       \cdot \left[t_{n-1} x^{n-1} y z 
       \sum_{k=\ldes(\pi)}^{\pi^{-1}(n-1)-1} (yz/x)^{k}\right.\\
& & + \left.x^{n-1} y z \sum_{k=\pi^{-1}(n-1)}^{n-2} 
      (yz/x)^{k} + y^{\ldes(\pi)} z^{n}\right]\\
&=& (1-yz/x)^{-1}\sum_{\pi\in T_{n-1}} t_{\Des(\pi^{-1})} x^{\inv(\pi)} \cdot\\
& & \cdot \left[t_{n-1} x^{n-1} y z 
    \left((yz/x)^{\ldes(\pi)} - (yz/x)^{\pi^{-1}(n-1)}\right)\right.\\
& & + \left. x^{n-1} y z \left((yz/x)^{\pi^{-1}(n-1)} - (yz/x)^{n-1}\right) 
    + y^{\ldes(\pi)} z^{n}(1-yz/x)\right]\\
&=& (1-yz/x)^{-1} [t_{n-1} x^{n-1} y z f_{n-1}(\bt,x,yz/x,1)\\
& & + \left.(1-t_{n-1}) x^{n-1} y z f_{n-1}(\bt,x,1,yz/x) - y^n z^n f_{n-1}(\bt,x,1,1)\right.\\
& &+ \left.(1-yz/x)z^n f_{n-1}(\bt,x,y,1)]\right..
\end{eqnarray*}
Multiply both sides by $(x-yz)$ to get the claimed recursion.

\qed

\begin{cor}\label{t.f_values}
The first few values of $f_n$ are:
\begin{eqnarray*}
f_1(\bt,x,y,z) &=& z\\
f_2(\bt,x,y,z) &=& z^2 + t_1xyz\\
f_3(\bt,x,y,z) &=& z^3 + t_1(x^2y^2z^2 + xyz^3) + t_2(x^2yz + xy^2z^2) \\
f_4(\bt,x,y,z) &=& z^4 + t_1(x^3y^3z^3 + x^2y^2z^4 + xyz^4)\\ 
               & & + \,t_2(x^4y^2z^2 + x^3y^3z^3 + x^2y^3z^3 + x^2yz^4 + xy^2z^4)\\
               & & + \,t_3(x^3yz + x^2y^2z^2 + xy^3z^3) 
                   + t_1t_3(x^3y^2z^2 + x^2y^3z^3)  
\end{eqnarray*}
\end{cor}


\section{Equidistribution of $\ldes$ and $\lind-1$}\label{s.equi}

In this section we prove Theorem~\ref{t.equi}.

\smallskip

\noindent
{\bf Note.}
Most of the permutations $\pi\in T_n$, namely those with $\pi(n)\ne n$,
satisfy $\lind(\pi) = \ldes(\pi)$. Nevertheless, the equidistributed parameters
are not $\ldes$ and $\lind$ but rather $\ldes$ and $\lind-1$.

\smallskip

\noindent
{\bf Note.}
$T_n(B)$ ``forgets'' whether or not $n-1$ belongs to $\Des(\pi^{-1})$.
The corresponding claim, without this ``forgetfulness'', is false!

\smallskip

\noindent
{\bf Proof of Theorem~\ref{t.equi}.}
We have to show that, letting $x=1$:
$$
q\hf_n(\bt,1,q,1) = \hf_n(\bt,1,1,q),
$$
where $\hf$ denotes $f$ under the additional substitution $t_{n-1}=1$.
This clearly holds for $n=1$ (as well as for $n=2,3,4$, by 
Corollary~\ref{t.f_values}).
By Theorem~\ref{t.rec}, for $n\ge 2$:
\begin{eqnarray*}
(1-yz) \hf_n(\bt,1,y,z) 
&=& yz f_{n-1}(\bt,1,yz,1) + (1-yz)z^n f_{n-1}(\bt,1,y,1)\\
& & - \,y^n z^n f_{n-1}(\bt,1,1,1).
\end{eqnarray*}
Letting $y=q$ and $z=1$ gives
\begin{eqnarray*}
(1-q) \hf_n(\bt,1,q,1) 
&=& q f_{n-1}(\bt,1,q,1) + (1-q) f_{n-1}(\bt,1,q,1)\\
& & - \,q^n f_{n-1}(\bt,1,1,1)\\
&=& f_{n-1}(\bt,1,q,1) - q^n f_{n-1}(\bt,1,1,1),
\end{eqnarray*}
whereas letting $y=1$ and $z=q$ yields
\begin{eqnarray*}
(1-q) \hf_n(\bt,1,1,q) 
&=& q f_{n-1}(\bt,1,q,1) + (1-q)q^n f_{n-1}(\bt,1,1,1)\\ 
& & - \,q^n f_{n-1}(\bt,1,1,1)\\ 
&=& q f_{n-1}(\bt,1,q,1) - q^{n+1} f_{n-1}(\bt,1,1,1).
\end{eqnarray*}
This completes the proof.
\qed

\section{Sign Balance of $T_n$}\label{s.balance}


In this section we prove Theorem~\ref{t.ldes.balance}.


\smallskip

\noindent
{\bf Proof of Theorem~\ref{t.ldes.balance}.}
Substituting $\bt = (1,1,\ldots)$, $z=1$, and $x=-1$ in 
the generating function~(\ref{e.gf}), denote
$$
g_n(-1,y) := f_n(\bar{1},-1,y,1) 
= \sum_{\pi\in T_n} (-1)^{\inv(\pi)} y^{\ldes(\pi)}\qquad(n\ge 1).
$$
This function records the sign and the last descent of $321$-avoiding permutations.
Recall also the generating function for last descent, from Theorem~\ref{t.equi}:
$$
g_n(1,y) := f_n(\bar{1},1,y,1) = \sum_{\pi\in T_n} y^{\ldes(\pi)}\qquad(n\ge 1).
$$
Let $g_0(1,y) := 1$.
We have to prove that
$$
g_{2n+1}(-1,y) = g_n(1,y^2)\qquad(n\ge 0),
$$
and
$$
g_{2n}(-1,y) = (1-y) g_n(1,y^2)\qquad(n\ge 1).
$$
By Theorem~\ref{t.rec}, the polynomial $g_n(-1,y)$ satisfies the recursion formula
\begin{eqnarray*}
(-1-y) g_n(-1,y) 
&=& (-1)^ny g_{n-1}(-1,-y) + (-1-y) g_{n-1}(-1,y)\\
& & + \,y^n g_{n-1}(-1,1).
\end{eqnarray*}
Clearly, 
$$
g_1(-1,y) = 1 = g_0(1,y^2)
$$
and also 
$$
g_2(-1,y) = 1-y = (1-y) g_1(1,y^2).
$$
We can proceed by induction.
If 
$$
g_{2n}(-1,y) = (1-y) g_n(1,y^2)
$$
then
\begin{eqnarray*}
(-1-y) g_{2n+1}(-1,y) 
&=& -y g_{2n}(-1,-y) + (-1-y) g_{2n}(-1,y)\\
& & + \,y^{2n+1} g_{2n}(-1,1)\\
&=& -y(1+y) g_n(1,y^2) + (-1-y)(1-y) g_n(1,y^2)\\
& & + \,0\\ 
&=& (-1-y) g_n(1,y^2),
\end{eqnarray*}
so that 
$$
g_{2n+1}(-1,y) = g_n(1,y^2).
$$
This, in turn, implies
\begin{eqnarray*}
(-1-y) g_{2n+2}(-1,y) 
&=& y g_{2n+1}(-1,-y) + (-1-y) g_{2n+1}(-1,y)\\
& & + \,y^{2n+2} g_{2n+1}(-1,1)\\
&=& y g_n(1,y^2) + (-1-y) g_n(1,y^2) + y^{2n+2} g_n(1,1)\\ 
&=& -g_n(1,y^2) + y^{2n+2} g_n(1,1) 
\end{eqnarray*}
and therefore
$$
g_{2n+2}(-1,y) = \frac{g_n(1,y^2) - y^{2n+2}g_n(1,1)}{1+y} = (1-y) h_n(y^2),
$$
where
$$
h_n(q) = \frac{g_n(1,q) - q^{n+1} g_n(1,1)}{1-q}.
$$
We need to show that 
$$
h_n(q) = g_{n+1}(1,q).
$$
Indeed, the equation
$$
(1-q) g_{n+1}(1,q) = g_n(1,q) - q^{n+1} g_n(1,1)
$$
follows immediately from Theorem~\ref{t.rec}.

\qed


\begin{cor}\label{t.balance} {\rm (equivalent to~\cite[Prop.~2]{SiS})}\\
The sign-balance enumerator for $T_n$
$$
g_n(-1,1) := f_n(\bar{1},-1,1,1) = \sum_{\pi\in T_n} (-1)^{\inv(\pi)}
$$
is either a Catalan number or zero:
\begin{eqnarray*}
g_{2n+1}(-1,1) &=& |T_n| = \frac{1}{n+1}{2n \choose n}\qquad(n\ge 0)\\
g_{2n}(-1,1) &=& 0 \qquad(n\ge 1)
\end{eqnarray*}
\end{cor}

The amazing connection between $T_n$, $T_{2n}$, and $T_{2n+1}$ 
calls for further study.

\section{Dyck Paths}\label{s.Dyck}

Let $\DP_n$ be the set of Dyck paths of length $2n$.
Thus, each $p\in \DP_n$ is a sequence $(p_1,\ldots,p_{2n})$ of
$n$ ``$+1$''s and $n$ ``$-1$''s, with nonnegative initial partial sums:
$$
p_1 + \ldots + p_i \ge 0\qquad(\forall i).
$$
Of course, $p_1 + \ldots + p_{2n} = 0$ by definition.

A {\em peak} of a Dyck path $p = (p_1,\ldots, p_{2n})\in \DP_n$ 
is an index $1\le i\le 2n-1$ such that $p_i = +1$ and $p_{i+1} = -1$.
Let $\Peak(p)$ be the set of all peaks of $p$.
Define the {\em Descent set} of $p$ to be
$$
\Des(p) := \{1\le i\le n-1 \,|\, i\in \Peak(p)\}.
$$
Define the {\em inverse path} $p^{-1}\in \DP_n$ to be
$$
p^{-1} := (-p_{2n},\ldots,-p_1),
$$
i.e., $p^{-1}$ is obtained by reversing the order of steps as well as the sign
of each step in $p$. Clearly,
$$
i\in \Peak (p^{-1}) \iff 2n-i\in \Peak(p),
$$
so that
\begin{eqnarray*}
\Des(p^{-1})
&=& \{1\le i\le n-1 \,|\, i\in \Peak(p^{-1})\}\\
&=& \{1\le i\le n-1 \,|\, 2n-i\in \Peak(p)\}.
\end{eqnarray*}
Thus the peaks strictly {\sl before} the midpoint of $p$ record the descents 
of $p$, whereas those strictly {\sl after} the midpoint 
of $p$ record the descents of $p^{-1}$, in reverse order. 
The reason for this terminology is the following result.
\begin{lem}\label{t.Dyck1}
There exists a bijection $\phi: T_n \to \DP_n$ such that,
if $\pi\in T_n$ and $p := \phi(\pi)\in\DP_n$, then
$$
\Des(p) = \Des(\pi)
$$
and
$$
\Des(p^{-1}) = \Des(\pi^{-1}). 
$$
\end{lem}

\noindent{\bf Proof.}
We shall define the bijection $\phi: T_n \to \DP_n$ in three steps.

%

First, let $\phi_1: T_n \to \SYT_2(n)$ be the Robinson-Schensted 
correspondence, where $\SYT_2(n)$ is the set of all pairs $(P,Q)$ of 
standard Young tableaux of size $n$ with the same two-rowed shape.

\smallskip
\noindent
{\bf Example.} 
$$
\phi_1(25134) =
\left(\matrix{1 & 3 & 4 & & 1 & 2 & 5 \cr 2 & 5 & & , & 3 & 4 & }\right) 
$$

Next, define $\phi_2: \SYT_2(n) \to \SYT(n,n)$, where $\SYT(n,n)$ is the set of
standard Young tableaux of (size $2n$ and) shape $(n,n)$.
$T=\phi_2(P,Q)$ is obtained by gluing $Q$ along its ``eastern frontier'' to 
the skew tableau obtained by rotating $P$ by $180^0$ and replacing each entry 
$1\le j\le n$ of $P$ by $2n+1-j$.

\smallskip
\noindent
{\bf Example.} 
$$
\phi_2 \left(\matrix{1 & 3 & 4 & & 1 & 2 & 5 \cr 2 & 5 & & , & 3 & 4 & }\right) =
\left.\matrix{1 & 2 & 5 & 6 & 9 \cr 3 & 4 & 7 & 8 & 10}\right.
$$

Finally, define $\phi_3: \SYT(n,n) \to \DP_n$ as follows:
the increasing ($+1$) steps in $\phi_3(T)$ are at places indexed by 
the entries in row $1$ of $T$, whereas the decreasing ($-1$) steps 
are at places indexed by row $2$ of $T$.

\smallskip
\noindent
{\bf Example.} 
$$
\phi_3 \left(\matrix{1 & 2 & 5 & 6 & 9 \cr 3 & 4 & 7 & 8 & 10}\right) =
(++--++--+-)
$$

\smallskip
The bijection we want is the composition $\phi = \phi_3 \phi_2 \phi_1$.

\medskip

Let us examine what happens to the descent sets of $\pi$ and of $\pi^{-1}$
under this sequence of transformations.

By a well-known property of the Robinson-Schensted correspondence $\phi_1$, if
$$
\pi \stackrel{\phi_1}{\longmapsto} 
(P,Q) \stackrel{\phi_2}{\longmapsto} 
T \stackrel{\phi_3}{\longmapsto} 
p 
$$
then 
$$
\pi^{-1} \stackrel{\phi_1}{\longmapsto}  
(Q,P) \stackrel{\phi_2}{\longmapsto} 
T^{-1} \stackrel{\phi_3}{\longmapsto} 
p^{-1}.
$$
Here $T^{-1}$ is the tableau obtained from $T\in \SYT(n,n)$ 
by a $180^0$ rotation and 
replacement of each entry $1\le j\le 2n$ by $2n+1-j$.

For 
a standard two-rowed Young tableau $P$ of size $n$ define:
$$
\Des(P) := \{1\le i\le n-1 \,|\,
\mbox{\rm $i$ is in row 1 and $i+1$ is in row 2 of $P$}\}. 
$$
For $T\in \SYT(n,n)$ (of size $2n$) define
\begin{eqnarray*}
\Des_1(T) &:=& \Des(T)\cap [n-1]\\
&=& \{1\le i\le n-1 \,|\, 
\mbox{\rm $i$ is in row 1 and $i+1$ is in row 2 of $T$}\}.
\end{eqnarray*}

Using all this notation we now have, for $\pi\in T_n$,
$$
\Des(\pi) = \Des(Q) = \Des_1(T) = \Des(p),
$$
where the equality $\Des(\pi) = \Des(Q)$ is again a well-known property of 
the Robinson-Schensted correspondence. Similarly,
$$
\Des(\pi^{-1}) = \Des(P) = \Des_1(T^{-1}) = \Des(p^{-1}).
$$
This completes the proof.

\qed

\smallskip

\noindent
{\bf Note.}
The bijection $\phi: T_n \to \DP_n$ in the proof of Lemma~\ref{t.Dyck1},
so well suited for our purposes, is essentially due to 
Knuth~\cite[p.~64]{Knuth} (up to an extra rotation of $T$).
Its ingredient $\phi_2$ was already used by MacMahon~\cite[p.~131]{MM}. 
For other bijections between these two sets see~\cite{BJS, Kratt, E, EP}.

\smallskip

For a Dyck path $p = (p_1,\ldots,p_{2n})\in \DP_n$, let
$$
\ldes(p):= \max \Des(p) = \max \{1\le i\le n-1 \,|\, p_i=+1, p_{i+1}=-1\}.
$$
If $\Des(p)=\emptyset$, i.e., if $p=(+\ldots+-\ldots-)$ is the ``unimodal path'',
we define $\ldes(p):=0$.

Define also
$$
\lind(p) := \max\{1\le i\le n \,|\, p_i=+1, p_{i+1}=-1\}.
$$
In other words,
$$
\lind(p) = \cases{n,& if $n\in \Peak(p)$;\cr
                  \ldes(p),& otherwise.}
$$                  

\begin{lem}\label{t.Dyck.ind}
Let
$$
\pi \stackrel{\phi_1}{\longmapsto} 
(P,Q) \stackrel{\phi_2}{\longmapsto} 
T \stackrel{\phi_3}{\longmapsto} 
p, 
$$
as in the proof of Lemma~\ref{t.Dyck1}.
The following conditions are equivalent:
\begin{enumerate}
\item
$\pi(n)=n$.
\item
$n$ is in row 1 of both $P$ and $Q$.
\item
$n$ is in row 1 of $T$, and $n+1$ is in row 2 of $T$. 
\item
$n\in \Peak(p)$.
\end{enumerate}
\end{lem}

\noindent{\bf Proof.}
If $\pi(n)=n$ then $n$ is the last element inserted into $P$, and it stays
in row 1 due to its size. Thus $n$ is in row 1 of both $P$ and $Q$.
Conversely, if $n$ is in row 1 of both $P$ and $Q$ then it is necessarily 
the {\em last} element in each of these rows. Thus it is the last element 
inserted into $P$, i.e., $\pi(n)=n$. 
This proves the equivalence of conditions 1 and 2.

Conditions 2 and 3 are clearly equivalent, by the definition of $\phi_2$.

Similarly for conditions 3 and 4, by the definitions of $\phi_3$ and of 
$\Peak(p)$.

\qed

\begin{lem}\label{t.Dyck3}
If $\pi\in T_n$ and $p := \phi(\pi)\in\DP_n$, then
$$
\ldes(p) = \ldes(\pi)
$$
and
$$
\lind(p) = \lind(\pi).
$$
\end{lem}

\noindent{\bf Proof.}
The first claim follows immediately from Lemma~\ref{t.Dyck1}, since
$$
\ldes(p) = \max \Des(p) = \max \Des(\pi) = \ldes(\pi).
$$
(If $\Des(p) = \Des(\pi) = \emptyset$ 
then $\ldes(p) = \ldes(\pi) = 0$ by definition.)

If $\pi(n)=n$ then $\lind(\pi)=\pi^{-1}(n)=n$, 
and also $\lind(p)=n$ by Lemma~\ref{t.Dyck.ind}.

Otherwise ($\pi(n)\ne n$) clearly $\pi^{-1}(n)\in \Des(\pi)$; 
and $\pi^{-1}(n)$ is actually the last descent of $\pi$, 
because of the $321$-avoiding condition. Thus $\ldes(\pi)=\lind(\pi)$.
On the other hand, $n\notin \Peak(p)$ by Lemma~\ref{t.Dyck.ind} so that, 
by definition, $\lind(p)=\ldes(p)$. Thus
$$
\lind(p) = \ldes(p) = \ldes(\pi) = \lind(\pi),
$$
as claimed.

\qed


\medskip
\noindent{\bf Proof of Theorem~\ref{t.Dyck2}.}
Combine Lemmas~\ref{t.Dyck1} and~\ref{t.Dyck3} with Theorem~\ref{t.equi}.

\qed

\section{Final Remarks}\label{remarks}

The quotient $P_n/ \langle QS_n \rangle$ of the polynomial ring 
$P_n=\bbq[x_1,\dots,x_n]$ by the ideal generated by
the quasi-symmetric functions without constant term
was studied by Aval, Bergeron and Bergeron~\cite{ABB}.
They determined its Hilbert series with respect to grading by total degree.
For a Dyck path $p\in \DP_n$, let
$$
\tail(p) := 2n - \max \{i \,|\, p_i=+1, p_{i+1}=-1\}
$$
be its ``tail length'', i.e., the number of consecutive decreasing steps at
its end. 
With this notation, the result of~\cite{ABB} can be formulated as follows.
\begin{thm}\label{t.ABB}{\rm \cite[Theorem 5.1]{ABB}}
The Hilbert series of the quotient $P_n/ \langle QS_n \rangle$  
(graded by total degree) is 
$$
\sum\limits_{p\in \DP_n} q^{n-\tail(p)} =
\sum_{k=0}^{n-1} \frac{n-k}{n+k} {n+k \choose k} q^k.
$$
\end{thm}


\begin{cor}\label{t.ldes1}
The Hilbert series of the quotient $P_n/ \langle QS_n \rangle$  
(graded by total degree) is equal to 
$$
\sum\limits_{\pi\in T_n} q^{\ldes(\pi)}. 
$$
\end{cor}

\noindent{\bf Proof.} 
By the definitions in Section~\ref{s.Dyck} above, 
for every Dyck path $p\in \DP_n$,
$$
\tail(p) = 2n - \max \Peak(p) = \min \Des(p^{-1}),
$$
with the convention $\min \emptyset := n$.
Thus, by Lemma~\ref{t.Dyck1}, if $\pi = \phi^{-1}(p)\in T_n$
then $\tail(p) = \min \Des(\pi^{-1})$. 
Now, define a bijection $\psi: T_n \to T_n$ by 
$$
\psi(\pi) := w_0 \pi^{-1} w_0\qquad(\forall \pi\in T_n),
$$
where $w_0 =n \ldots 21$ is the longest element in $S_n$.
If 
$$
\pi \stackrel{\psi}{\longmapsto} 
\hpi \stackrel{\phi}{\longmapsto} p
$$
then clearly
$$
\tail(p) = \min \Des(\hpi^{-1}) = n - \max \Des(\pi) = n - \ldes(\pi).
$$
We conclude that
$$
\sum\limits_{p\in \DP_n} q^{n-\tail(p)} =
\sum\limits_{\pi\in T_n} q^{\ldes(\pi)}.
$$

\qed


%

\medskip
Theorem~\ref{t.ABB} was proved in~\cite{ABB} by constructing a monomial basis 
for the quotient ring $P_n/ \langle QS_n \rangle$. By the proof of 
Corollary~\ref{t.ldes1}, the elements of this basis may be indexed by 
permutations $\pi\in T_n$, with total degree given by $\ldes(\pi)$.
It is well known that the Hilbert series of
the graded 
coinvariant algebra $P_n/I_n$, where $I_n$ is the ideal generated by the 
symmetric functions without constant term, is equal to
$$
\sum_{\pi\in S_n} q^{\maj(\pi)},
$$
where $\maj(\pi):=\sum_{i\in \Des(\pi)} i$.
A monomial basis idexed by permutations $\pi\in S_n$, with total degree 
given by $\maj(\pi)$, was constructed in~\cite{GS}.
In this sense, the last descent may be considered as a $T_n$-analogue 
of the major index.

\begin{prb}
Find algebraic interpretations and proofs for Theorems~\ref{t.equi} and~\ref{t.ldes.balance}.
\end{prb}

For a different algebraic interpretation of the generating function
in Theorem~\ref{t.ABB} see \cite{Panyushev}.


\bigskip

\noindent{\bf Acknowledgments.} Thanks to Nantel Bergeron for stimulating 
discussions, and to Astrid Reifegerste and an anonymous referee for their comments.

\end{document}